\documentclass[12pt,a4paper]{article}
\usepackage{amsmath,amssymb,amsfonts}
\def\Bbb{\mathbb}

\title{\bf On a criterion for the equality of Dedekind sums}

\author{Kurt Girstmair}
\date{}

\makeatletter
\let\@@maketitle=\maketitle
\def\maketitle{\def\thispagestyle##1{\relax}\@@maketitle}
\makeatother
\newtheorem{theorem}{Theorem}

\def\BE{\begin{equation}}
\def\EE{\end{equation}}
\def\BD{\begin{displaymath}}
\def\ED{\end{displaymath}}
\def\BA{\begin{array}}
\def\EA{\end{array}}
\def\BEA{\begin{eqnarray*}}
\def\EEA{\end{eqnarray*}}
\def\BI{\bibitem}

\def\Z{\Bbb Z}

\def\phi{\varphi}

\def\CMOD#1#2#3{#1 \equiv #2 \: \mbox{mod}\: #3}

\def\MB{\mbox}
\def\LD{\ldots}

\def\MN{\medskip\noindent}
\def\STOP{\hfill$\Box$}

\def\DED{Dedekind }

\begin{document}
\maketitle

\begin{abstract}
\noindent
In \cite{Ja} it was shown that the \DED sums $s(m_1,n)$ and $s(m_2,n)$ are equal only if
$(m_1m_2-1)(m_1-m_2)\equiv 0$ mod $n$. Here we show that the latter condition is equivalent to
$12s(m_1,n)-12s(m_2,n)\in \Z$. In addition, we determine, for a given number $m_1$,
the number of integers $m_2$ in the range $0\le m_2<n$, $(m_1,m_2)=1$, such that $12s(m_1,n)-12s(m_2,n)\in \Z$,
provided that $n$ is square-free.
\end{abstract}

\section*{1. Introduction and results}

Let $n$ be a positive integer and $m\in \Z$, $(m,n)=1$. The classical \DED sum $s(m/n)$ is defined by
\BD
   s(m/n)=\sum_{k=1}^n ((k/n))((mk/n))
\ED
where $((\LD))$ is the usual sawtooth function (see, for instance, \cite{RaGr}, p. 1). In the present setting it is more
natural to work with
\BD
S(m/n)=12s(m/n)
\ED instead. In \cite{Ja} the following theorem was shown:

\begin{theorem} 
\label{t1}
Let $m_1, m_2$ be integers that are relatively prime to $n$. If $S(m_1/n)=S(m_2/n)$,
then
\BE
\label{2}
(m_1m_2-1)(m_1-m_2)\equiv 0 \MB{ mod }n.
\EE
\end{theorem} 

\MN
So (\ref{2}) is only a {\em necessary} condition for the equality of $S(m_1,n)$
and $S(m_2,n)$. But what does this condition really stand for? In this note we show

\begin{theorem} 
\label{t2}
Let $m_1, m_2$ be integers that are relatively prime to $n$.
Then $S(m_1/n)-S(m_2/n)\in \Z$ if, and only if,
{\rm (\ref{2})} holds.
\end{theorem} 

Suppose that the number $m_1$, $(m_1,n)=1$, is given. Observe that $S(m_2,n)$ depends only on
the residue class of $m_2$ mod $n$. Hence it is natural to ask how many numbers $m_2$, $(m_2,n)=1$, exist in the
range $0\le m_2<n$ such that $S(m_1,n)$ and $S(m_2,n)$ have equal fractional parts. We consider only a
simple case here, namely,

\begin{theorem} 
\label{t3}
Let $n=p_1\LD p_t$ be square-free, so $p_1,\LD,p_t$ are pairwise different primes. For a given
number $m_1$, $(m_1,n)=1$, we have
\BD
|\{m_2\,:\,0\le m_2<n,\,(m_2,n)=1,\, S(m_1/n)-S(m_2/n)\in \Z\}|=2^s,
\ED
where $s=|\{j\,:\,1\le j\le t,\,m_1\not\equiv \pm 1 \MB{ mod } p_j\}|$.
\end{theorem} 

\MN
{\em Example.} Let $n=15015=3\cdot 5\cdot 7\cdot 11\cdot 13$ and $m_1=17$.
Clearly, $m_1\equiv -1$ mod 3 but $m_1\not\equiv\pm 1$ mod $p$ for $p\in\{5,7,11,13\}$.
So Theorem \ref{t3} says that there are $2^4=16$ numbers $m_2$, $0\le m_2<n$, $(m_2,n)=1$,
such that $S(m_2,n)$ has the same fractional part as $S(17,n)$. In fact, we obtain
\BD
  S(m_2,n)=\frac{710}{3003}+\left\{\begin{array}{rl}
               880, & \MB{ for }m_2\in\{17, 3533\}, \\
               40,  & \MB{ for }m_2\in\{6023,12542\},  \\
               16,  & \MB{ for }m_2\in\{992,2558,6452,12113\}, \\
               -8,  & \MB{ for }m_2\in\{563,2987,6107,6998,11567,12458\},  \\
               -32, & \MB{ for }m_2\in\{8993,9572\}.
             \end{array}\right.
\ED

\section*{2. Proofs}

{\em Proof of Theorem \ref{t2}.}
Let $m$ be an integer, $(m,n)=1$. In a first step we use the Barkan-Hickerson-Knuth formula in order
to determine the fractional part of $S(m,n)$.
We start with the continued fraction expansion $m/n=[a_0,a_1,\LD,a_k]$. Since $S(m/n)$ depends only on
the residue class of $m$ mod $n$,
we may assume $0\le m<n$, i. e., $a_0=0.$
Let $s_0/t_0, \LD ,s_k/t_k=m/n$ be the convergents of $m/n$, where the numbers $s_j$, $t_j$, $0\le j\le k$,
are recursively defined as usual (see \cite{RoSz}, p. 2). In particular, $s_k=m$ and $t_k=n$. The
Barkan-Hickerson-Knuth formula says that for $k\ge 1$
\BD
\label{3}
S(s_k/t_k)=\sum_{j=1}^k(-1)^{j-1}a_j+\left\{\begin{array}{ll} \vspace{5mm}
                                                              (s_k+ t_{k-1})/{t_k}-3 & \MB{if } k \MB{ is odd}, \\
                                                               (s_k-t_{k-1})/{t_k}  & \MB{if } k \MB{ is even,}
                                                            \end{array}\right.
\ED
see \cite{Ba}, \cite{Hi}, \cite{Kn}.
Hence,
\BE
\label{4}
S(m,n)\equiv \frac{m+(-1)^{k-1}t_{k-1}}n \MB{ mod }\Z.
\EE
Further, we observe the basic identity
\BD
\label{5}
  s_kt_{k-1}-t_ks_{k-1}=(-1)^{k-1},
\ED
which gives
\BE
\label{6}
  \CMOD{mt_{k-1}}{(-1)^{k-1}}n
\EE
(see \cite{RoSz}, p. 2). Let $m^*$ denote the {\em inverse} of $m$ mod $n$, i. e., the uniquely determined
integer $m^*$ in the range $0\le m^*<n$ such that $mm^*\equiv 1$ mod $n$.
Then (\ref{6}) means
\BD
   \frac{(-1)^{k-1}t_{k-1}}n\equiv \frac{m^*}n \MB{ mod }\Z.
\ED
Together with (\ref{4}), this gives
\BE
\label{8}
S(m,n)\equiv \frac{m+m^*}n \MB{ mod }\Z.
\EE
This congruence is also valid in the (trivial) case $k=0$, where $n=1$ and $S(m,n)=0$.

The second (and final) step of the proof of Theorem \ref{t2} consists in showing
that (\ref{2}) is equivalent to $S(m_1,n)\equiv S(m_2,n)$ mod $\Z$.
To this end we note that both the \DED sums $S(m_1,n)$, $S(m_2,n)$ and the condition (\ref{2})
depend only on the residue classes of $m_1$ mod $n$ and $m_2$ mod $n$. Hence we may assume
$0\le m_1,m_2<n$. By (\ref{8}),
$S(m_1/n)\equiv S(m_2,n)$ mod $\Z$ if, and only if,
\BE
\label{12}
 m_1+m_1^*\equiv m_2+m_2^* \MB{ mod } n.
\EE
The proof is complete if we can show that (\ref{12}) is equivalent to (\ref{2}).
However, multiplying (\ref{12}) by $m_1m_2$, we obtain
\BE
\label{14}
m_1^2m_2+m_2\equiv m_1m_2^2+m_1 \MB{ mod n }
\EE
and
\BD
m_1m_2(m_1-m_2)\equiv m_1-m_2 \MB{ mod n},
\ED
which obviously yields (\ref{2}). Conversely, (\ref{2}) implies (\ref{14}),
and on multiplying this congruence by $m_1^*m_2^*$, we obtain (\ref{12}).
\STOP

\MN
{\em Proof of Theorem 3.}
By Theorem \ref{t2},
\BD
\{m_2\,:\,0\le m_2<n,\,(m_2,n)=1,\, S(m_1/n)-S(m_2/n)\in \Z\}=
\ED \BD
\{m_2\,:\,0\le m_2<n,\,(m_2,n)=1,\,(m_1m_2-1)(m_1-m_2)\equiv 0 \MB{ mod }n\}.
\ED
By the Chinese Remainder Theorem, the latter set has the cardinality
\BE
\label{16}
   \prod_{j=1}^tL_j,
\EE
where $L_j=|\{m_2\,:\,1\le m_2<p_j,\,(m_1m_2-1)(m_1-m_2)\equiv 0 \MB{ mod }p_j\}|$,
$j=1,\LD,t$. Now it is easy to see that
\BD
   L_j=\left\{\begin{array}{rl}
                1 & \MB{ if } m_1\equiv \pm 1 \MB{ mod }p_j, \\
                2 &  \MB{ otherwise.}
              \end{array}\right.
\ED
Therefore (\ref{16}) gives the desired result.
\STOP


\vspace{0.5cm}
\noindent
Kurt Girstmair            \\
Institut f\"ur Mathematik \\
Universit\"at Innsbruck   \\
Technikerstr. 13/7        \\
A-6020 Innsbruck, Austria \\
Kurt.Girstmair@uibk.ac.at

\end{document}